\documentclass [11pt,a4paper]{article}
\usepackage{amssymb}
\def\F{\mathbb{F}}

\def\Z{\mathbb{Z}}

\newtheorem{theorem}{Theorem}
\newtheorem{lemma}{Lemma}

\newcommand{\quash}[1]{}

\begin{document}

\title{Linear complexity problems of level sequences of Euler quotients and their related binary sequences}
\author{
Zhihua Niu\\
School of Computer Engineering and Science, Shanghai University,\\
Shangda Road, Shanghai 200444, P. R. China \\
zhniu@staff.shu.edu.cn\\
\\
Zhixiong Chen\\
School of Mathematics, Putian University, \\ Putian, Fujian
351100, P. R. China\\
ptczx@126.com\\
\\
Xiaoni Du\\
 College of Mathematics and Information Science, \\ Northwest Normal University,
Lanzhou, Gansu 730070, P. R. China
 }

\maketitle

\begin{abstract}
The Euler quotient modulo an odd-prime power $p^r~(r>1)$ can be uniquely
decomposed as a $p$-adic number of the form
$$
\frac{u^{(p-1)p^{r-1}} -1}{p^r}\equiv a_0(u)+a_1(u)p+\ldots+a_{r-1}(u)p^{r-1} \pmod {p^r},~ \gcd(u,p)=1,
$$
where $0\le a_j(u)<p$ for $0\le j\le r-1$  and  we set all $a_j(u)=0$ if $\gcd(u,p)>1$. We firstly study certain arithmetic properties of the level sequences $(a_j(u))_{u\ge 0}$ over $\F_p$ via introducing a new quotient. Then we determine the exact values of linear complexity of $(a_j(u))_{u\ge 0}$ and values of $k$-error linear complexity
for binary sequences defined by $(a_j(u))_{u\ge 0}$.
\end{abstract}



\section{Introduction}

Let $p$ be an odd prime and $r$ be a positive integer. For all integers $u$ with $\gcd(u,p)=1$, by the Euler Theorem  we have
$$
u^{\varphi(p^r)}\equiv 1 \pmod {p^r},
$$
where $\varphi(-)$ is the Euler Totient function. Hence
we define $Q_{r}(u)$  modulo $p^r$   by
\begin{equation}\label{def}
 Q_{r}(u) \equiv \frac{u^{\varphi(p^r)} -1}{p^r} \pmod {p^r}, ~0 \le
Q_{r}(u) < p^r,~\mathrm{if}~\gcd(u,p)=1,
\end{equation}
which is called the \emph{Euler quotient} in \cite{ADS}. In fact, if we write
\begin{equation}\label{def-2}
u^{\varphi(p^r)}= 1+a_1p^r+a_2p^{2r}+\ldots\in\mathbb{Z}, ~ 0\le a_i<p^r ~ \mathrm{for}~ i\ge  1,
\end{equation}
we have $Q_{r}(u)=a_1$.
For convenience, we set
\begin{equation}\label{zero}
Q_{r}(lp) = 0, \qquad l \in \mathbb{Z}.
\end{equation}
If $r=1$, $Q_{1}(u)$ is also called the \emph{Fermat quotient}. A more general notion, called the \emph{Carmichael  Quotient},
is studied in \cite{Sha}.
Many number theoretic questions have been studied for these quotients
and their generalizations
\cite{ADS,BFKS,C,COW,CW,CW2,CW3,EM,GW,OS,Sha,Shk,S,S2010,S2011,S2011b,SW}.

Let $\mathbb{Z}_{p^r}$ be the integer residue ring modulo $p^r$. Any element $a\in\mathbb{Z}_{p^r}$ has a unique $p$-adic decomposition as $a = a_0 + a_1  p + \ldots +a_{r-1}p^{r-1}$, where $a_i\in\{0, 1, \ldots, p-1\}$. Hence for a sequence $(s(u))_{u\ge 0}$ over $\mathbb{Z}_{p^r}$, it has a unique
$p$-adic decomposition as
$$
s(u) = s_0(u) + s_1(u)  p + \ldots +s_{r-1}(u)p^{r-1},~ u\ge 0,
$$
where $(s_i(u))_{u\ge 0}$ is a sequence over $\{0, 1, \ldots, p-1\}$. The sequence $(s_i(u))_{u\ge 0}$ is called the \emph{$i$-th level sequence} of $(s(u))_{u\ge 0}$,
and $(s_{r-1}(u))_{u\ge 0}$ the \emph{highest-level sequence} of $(s(u))_{u\ge 0}$. They can be naturally considered as
the sequences over the finite field $\F_p$. S. Q. Fan and W. F. Qi (partly with coauthors) extensively investigated the level sequences of linear recurring sequences over
$\mathbb{Z}_{p^r}$ (or more generally $\mathbb{Z}_{M}$, where $M>1$ is an arbitrary number),  see \cite{FH03,FH,ZQT,ZQ10,ZQ05,ZQ07,QZ} and references therein. Certain $(s(u))_{u\ge 0}$ over $\mathbb{Z}_{p^r}$ is relevant to
FCSR sequences \cite{TQ}.

On the other hand,  Fermat quotients, Euler quotients and  Carmichael  Quotients have been studied recently from the
viewpoint of cryptography, see \cite{AW,CW4,CD,CDR2014,CG,CHD,CN,COW,DCH,DKC,GW,OS,WCD-12}. More exactly, the authors of \cite{OS} studied the linear complexity profile of the Fermat quotient sequence $(Q_{1}(u))_{u\ge 0}$. As we know, this is the first work to consider the cryptographic feature of Fermat quotients. The authors of \cite{COW,GW} used Fermat quotients and Euler quotients to define pseudorandom sequences.
The first one is the binary threshold sequence $(e(u))_{u\ge 0}$  defined by
\begin{equation}\label{binarythreshold}
e(u)=\left\{
\begin{array}{ll}
0, & \mathrm{if}\,\ 0\leq Q_{r}(u)/p^r< \frac{1}{2},\\
1, & \mathrm{if}\,\ \frac{1}{2}\leq Q_{r}(u)/p^r< 1,
\end{array}
\right. \quad  u \ge 0.
\end{equation}
The second one, by combining $Q_{r}(u)$ with $\chi$, which is a fixed multiplicative
character modulo $p^r$ of order $m>1$, is the $m$-ary sequences $(\widetilde{e}(u))_{u\ge 0}$  defined by
\begin{equation}\label{m-ary}
  \exp(2\pi i \widetilde{e}(u)/m)=
     \chi(Q_{r}(u)),~ 0\le \widetilde{e}(u)<m \quad \mbox{if } \gcd(Q_{r}(u), p)=1
\end{equation}
and $\widetilde{e}(u)=0$ otherwise. Most recent studies are concentrated in the case of $r=1$: the authors of \cite{COW,GW} investigated measures of pseudorandomness as well as linear complexity profile of
$(e(u))_{u\ge 0}$ and $(\widetilde{e}(u))_{u\ge 0}$ via certain character sums over Fermat quotients,
The authors of \cite{CHD,DKC}  determined the \emph{linear complexity} (see
 below for the definition) of $(e(u))_{u\ge 0}$ and $(\widetilde{e}(u))_{u\ge 0}$ if $2$ is a primitive element modulo $p^2$, and later the authors of \cite{CD,CG,CW4}  extended to a more general setting of $2^{p-1}\not\equiv 1 \pmod {p^2}$, the authors of \cite{CW4,CN}  also determined the  trace representations and the \emph{$k$-error linear complexity} (see  below for the definition) of $(e(u))_{u\ge 0}$ and $(\widetilde{e}(u))_{u\ge 0}$, respectively. The authors of \cite{DCH} extended \cite{CD} furtherly to determine the linear complexity  of $(e(u))_{u\ge 0}$ when $r>1$  under the assumption of $2^{p-1}\not\equiv 1 \pmod {p^2}$.  We refer the reader to related references for details. All results indicate that such sequences have desirable cryptographic features.

Hence in this paper, we describe the Euler quotient $Q_{r}(u)$ as the $p$-adic decomposition
$$
Q_{r}(u)=a_0(u)+a_1(u)p+\ldots+a_{r-1}(u)p^{r-1}, ~ u\ge 0
$$
where $0\le a_j(u)<p$ for $0\le j\le r-1$, and consider the linear complexity of the level sequences $(a_j(u))_{u\ge 0}$ over $\F_p$ via introducing a new quotient, which coincides with the level sequences $(a_j(u))_{u\ge 0}$.  Our second aim is to determine the $k$-error linear complexity
for certain binary sequences defined by the level sequences $(a_j(u))_{u\ge 0}$ of the Euler quotient $Q_{r}(u)$.

We conclude this section by recalling the notions of the linear
complexity and the $k$-error linear complexity. Let $\F$ be a field.  For a $T$-periodic
sequence $(s(u))_{u\ge 0}$ over $\F$, we recall that the
\emph{linear complexity} over $\F$, denoted by  $LC^{\F}((s(u))_{u\ge 0})$, is the least order $L$ of a linear
recurrence relation over $\mathbb{F}$
$$
s(u+L) = c_{L-1}s(u+L-1) +\ldots +c_1s(u+1)+ c_0s(u)\quad
\mathrm{for}\,\ u \geq 0,
$$
which is satisfied by $(s(u))_{u\ge 0}$ and where $c_0\neq 0, c_1, \ldots,
c_{L-1}\in \mathbb{F}$.
Let
$$
S(X)=s(0)+s(1)X+s(2)X^2+\ldots+s(T-1)X^{T-1}\in \mathbb{F}[X],
$$
which is called the \emph{generating polynomial} of $(s(u))_{u\ge 0}$. Then the linear
complexity over $\F$ of $(s(u))_{u\ge 0}$ is computed by
\begin{equation}\label{licom}
  LC^{\F}((s(u))_{u\ge 0}) =T-\deg\left(\mathrm{gcd}(X^T-1,
  ~S(X))\right),
\end{equation}
see, e.g. \cite{CDR,LN} for details. For integers $k\ge 0$, the \emph{$k$-error linear complexity} over $\F$ of $(s(u))_{u\ge 0}$, denoted by $LC^{\F}_k((s(u))_{u\ge 0})$, is the smallest linear complexity (over $\F$) that can be
obtained by changing at most $k$ terms of the sequence per period, see \cite{SM,Meidl}, and see \cite{DXS} for the related even earlier defined \emph{sphere complexity}.  Clearly $LC^{\F}_0((s(u))_{u\ge 0})=LC^{\F}((s(u))_{u\ge 0})$ and
$$
T\ge LC^{\F}_0((s(u))_{u\ge 0})\ge LC^{\F}_1((s(u))_{u\ge 0})\ge \ldots \ge LC^{\F}_k((s(u))_{u\ge 0})=0
$$
when $k$ equals the number of nonzero terms of $(s(u))_{u\ge 0}$ per period, i.e., the weight of $(s(u))_{u\ge 0}$.

The linear complexity and the $k$-error linear complexity are important cryptographic characteristics of sequences
and provide information on the predictability and thus unsuitability for cryptography. For a sequence to be cryptographically strong, its linear complexity
should be large, but not significantly reduced by changing a few
terms. And according to  the Berlekamp-Massey
algorithm \cite{Massey}, the linear complexity
should be at least a half of the period.

\section{A new quotient}

In this section, we introduce a new quotient to represent the level sequences of  the Euler quotient $Q_{r}(u)$.

For integers $r> s>0$, we can check
\begin{equation}\label{q-r-s}
Q_{r}(u)\equiv Q_{s}(u) \bmod {p^s}, ~~ u\ge 0.
\end{equation}
In fact, for $p|u$ we have $Q_{r}(u)= Q_{s}(u)=0$ by the assumption of (\ref{zero}). Now we suppose $\gcd(u,p)=1$.
Let
$$
u^{\varphi(p^s)}=1+b_1p^s+b_2p^{2s}+\ldots \in \mathbb{Z}, ~~ 0\le b_1,b_2,\ldots< p^s,
$$
we see that $Q_{s}(u)=b_1$ by (\ref{def-2}). On the other hand, we verify
\begin{eqnarray*}
u^{\varphi(p^r)}& = & (u^{\varphi(p^s)})^{p^{r-s}}=
(1+b_1p^s+b_2p^{2s}+\ldots)^{p^{r-s}}\\
& = & 1+b_1p^r+(b_1^2 b_2(p^{r-s}-1)/2)p^{r+s}+\ldots \in \mathbb{Z},
\end{eqnarray*}
from which we derive
$$
\frac{u^{\varphi(p^r)} -1}{p^r}=b_1+(b_1^2 b_2(p^{r-s}-1)/2)p^{s}+\ldots.
$$
We get $Q_{r}(u)\equiv b_1 \bmod {p^s}$. Hence we prove (\ref{q-r-s}).

 From (\ref{q-r-s}), for integer $r\ge 2$ one can define a new quotient by
\begin{equation}\label{new-quotient}
H_{r-1}(u)\equiv \frac{Q_{r}(u)-Q_{r-1}(u)}{p^{r-1}}\bmod {p}, ~~ 0\le H_{r-1}(u)< p, ~~ u\ge 0.
\end{equation}
Indeed, we can write
$$
Q_{r}(u)=Q_{r-1}(u)+H_{r-1}(u)p^{r-1}, ~ u\ge 0
$$
and hence
$$
Q_{r}(u)=H_{0}(u)+H_{1}(u)p+\ldots +H_{r-1}(u)p^{r-1}, ~ u\ge 0
$$
by induction on $r-1$, where $H_{0}(u)=Q_{1}(u)$ by (\ref{q-r-s}). Hence $(H_{i}(u))_{u\ge 0}$ is indeed the highest level
sequence of $(Q_{i+1}(u))_{u\ge 0}$ for $i\ge 1$.

For example, if $r=2$ and
$$
u^{p-1}=1+c_1p+c_2p^{2}+\cdots \in \mathbb{Z}, ~~ 0\le c_1,c_2,\ldots< p,
$$
we have $Q_{1}(u)=c_1$ and $Q_{2}(u)=c_1+(\frac{p-1}{2}c_1^2+c_2)p \pmod {p^2}$, and hence
$$
H_{0}(u)=c_1, ~~ H_{1}(u)\equiv \frac{p-1}{2}c_1^2+c_2 \pmod {p}.
$$
For $p|u$, we have $H_{0}(u)=H_{1}(u)=0$.

Since the $i$-th level sequence of $(Q_{r}(u))_{u\ge 0}$ is the highest level
sequence of $(Q_{i+1}(u))_{u\ge 0}$ for $i\ge 0$, we only consider the highest level
sequence of $(Q_{r}(u))_{u\ge 0}$ in the context, i.e., the quotient $H_{r-1}(u)$. Below we prove two simple properties for
$H_{r-1}(u)$. We remark again that $H_{0}(u)=Q_{1}(u)$, which is the Fermat quotient.

\begin{theorem}\label{HH}
For any integers $v$, $k$ and $r\ge 1$, we have
$$
H_{r-1}(v+kp^{r})\equiv H_{r-1}(v) - kv^{p-2} \bmod {p}.
$$
\end{theorem}
Proof. For $r=1$, $H_{0}(u)$ is the Fermat quotient $Q_{1}(u)$ and the result follows, see  \cite{OS}.
For $r>1$, since the least period of $(Q_{r-1}(u))_{u\ge 0}$ is $p^{r}$, together with $v^p\equiv v \bmod p$ we get
\begin{eqnarray*}
H_{r-1}(v+kp^{r}) & \equiv & \frac{Q_{r}(v+kp^{r})-Q_{r-1}(v)}{p^{r-1}}\\
                    & \equiv & \frac{Q_{r}(v)-Q_{r-1}(v)}{p^{r-1}}+k(p-1)v^{\varphi(p^r)-1}\\
                    & \equiv & H_{r-1}(v)+k(p-1)v^{p-2}\bmod {p}.
\end{eqnarray*}
We complete the proof. ~\hfill $\Box$

The least period of $(H_{r-1}(u))_{u\ge 0}$ follows from Theorem \ref{HH} directly.

\begin{theorem}\label{HH-period}
For integer $r\ge 1$,  the least period of $(H_{r-1}(u))_{u\ge 0}$  is $p^{r+1}$.
\end{theorem}

We remark that W. Leeb \cite{L} extended the Fermat quotients to introduce the notion of  \emph{Fermat quotients of order $i\ge 1$} by defining
$$
F^{(1)}(u)= Q_{1}(u)
$$
and for $i>1$
\begin{equation}\label{highorder}
F^{(i)}(u)\equiv \frac{u^{p-1}-1-F^{(1)}(u)p-F^{(2)}(u)p^2-\ldots -F^{(i-1)}(u)p^{i-1}}{p^i} \bmod p,
\end{equation}
with $0\le F^{(i)}(u)<p$ for all integers $u$ with $\gcd(u,p)=1$ and $F^{(i)}(u)=0$ otherwise. Indeed,
$$
F^{(i)}(u)=c_i, ~~ i\ge 1
$$
for $\gcd(u,p)=1$, if
$$
u^{p-1}=1+c_1p+c_2p^{2}+\cdots \in \mathbb{Z}, ~~ 0\le c_1,c_2,\ldots< p.
$$
We find that $F^{(i)}(u)$ is different from $H_{r-1}(u)$ defined in (\ref{new-quotient}).
 (Note that W. Leeb introduced this definition for more general settings.)

\section{Linear complexity of level sequences}

In this section, we determine the exact value of the linear complexity of the highest-level sequence $(H_{r-1}(u))_{u\ge 0}$ of the Euler quotient $Q_{r}(u)$.

\begin{theorem}\label{LC-p}
For integers $r\ge 1$,  the linear complexity (over the finite field $\F_p$) of the highest-level sequence $(H_{r-1}(u))_{u\ge 0}$ of Euler quotients in (\ref{def}) and (\ref{zero}) satisfies
$$
LC^{\F_p}((H_{r-1}(u))_{u\ge 0})=p^{r}+p-1.
$$
\end{theorem}
Proof. From Theorem \ref{HH-period}, the least period of $(H_{r-1}(u))_{u\ge 0}$  is $p^{r+1}$. So for all integers $u\equiv i_0+i_1p+\ldots+i_{r}p^{r} \pmod {p^{r+1}}$ with $0\le i_0, i_1, \ldots, i_{r}<p$,  we see that $(H_{r-1}(u))_{u\ge 0}$ can be represented by
$$
H_{r-1}(i_0+i_1p+\ldots+i_{r}p^{r}+jp^{r+1})=\rho(i_0, i_1, \ldots, i_{r}) ~~\mathrm{for}~~ j\ge 0,
$$
where the polynomial $\rho(X_0,X_1,\ldots, X_{r})\in \F_p[X_0,X_1,\ldots, X_{r}]/\langle X^p_0-X_0, X^p_1-X_1, \ldots, X^p_r-X_{r}\rangle$ is of the form
\begin{eqnarray*}
&& \rho(X_0,X_1,\ldots, X_{r}) \\
& = & \sum\limits_{c_0=0}^{p-1}\sum\limits_{c_1=0}^{p-1}\cdots\sum\limits_{c_{r-1}=0}^{p-1}
H_{r-1}(c_0+c_1p+\ldots +c_{r-1}p^{r-1})\prod_{l=0}^{r-1}\left(1-(X_l-c_l)^{p-1} \right) - X_{r}X_0^{p-2},
\end{eqnarray*}
since
\begin{eqnarray*}
& &H_{r-1}(i_0+i_1p+\ldots+i_{r}p^{r}+jp^{r+1})\\
&\equiv & H_{r-1}(i_0+i_1p+\ldots+i_{r-1}p^{r-1}) - i_{r}(i_0+i_1p+\ldots+i_{r-1}p^{r-1})^{p-2}\\
&\equiv & \sum\limits_{c_0=0}^{p-1}\sum\limits_{c_1=0}^{p-1}\cdots\sum\limits_{c_{r-1}=0}^{p-1}
H_{r-1}(c_0+c_1p+\cdots +c_{r-1}p^{r-1})\prod_{l=0}^{r-1}\left(1-(i_l-c_l)^{p-1} \right) - i_{r}i_0^{p-2} \bmod {p}
\end{eqnarray*}
by Theorem \ref{HH}.

 Then the degree of $\rho(X_0,X_1,\ldots, X_{r})$ is $\deg(\rho)=p^{r}+p-2$, see \cite{BEP} for the definition of the degree of multi-variable polynomials. Hence by \cite[Theorem 8]{BEP}, we have $LC^{\F_p}((H_{r-1}(u))_{u\ge 0})=\deg(\rho)+1=p^{r}+p-1$.
 ~\hfill $\Box$

 The case of $r=1$ in Theorem \ref{LC-p} has been reported in \cite{OS}.

\section{Linear complexity and $k$-error linear complexity of binary sequences derived from level sequences}

In this section, we apply the highest-level sequence $(H_{r-1}(u))_{u\ge 0}$ of the Euler quotient $Q_{r}(u)$ to defining
some families of binary sequences and determine their linear complexity and $k$-error linear complexity. Suppose that $2$ is a primitive root modulo $p^2$. Then it is clear that $2$ is also a primitive root modulo $p^n$ for every $n\ge 1$, see e.g. \cite{Na}.

From Theorem \ref{HH}, the quotient
$H_{r-1}(-)$ induces a surjective map from $\Z_{p^{r+1}}^*$ (the group of invertible
elements modulo $p^{r+1}$) to $\Z_p$ (the additive group of numbers modulo $p$). Let
$$
D_l=\{u: 0\le u< p^{r+1},~ \gcd(u,p)=1,~H_{r-1}(u)=l\}
$$
for $l=0,1,\ldots,p-1$ and $P=\{lp : 0\le l< p^r\}$. We define a $p^{r+1}$-periodic binary sequence $(f(u))_{u\ge 0}$ by
\begin{equation}\label{f-binary}
f(u)=\left\{
\begin{array}{ll}
1, & \mathrm{if}\,\ u\bmod {p^{r+1}}\in \cup_{l\in \mathcal{I}} D_l,\\
0, & \mathrm{otherwise},
\end{array}
\right. \quad u\ge 0,
\end{equation}
 where $\mathcal{I}$ is a non-empty subset of $\{0,1,\ldots,p-1\}$.  In particular, if $\mathcal{I}=\{\frac{p+1}{2}, \frac{p+1}{2}+1, \ldots, p-1\}$, $(f(u))_{u\ge 0}$ is the binary threshold sequence defined in (\ref{binarythreshold}) when $r=1$ and
if $\mathcal{I}$ is the set of quadratic non-residues modulo $p$, $(f(u))_{u\ge 0}$ is the binary sequence defined in (\ref{m-ary}) when  $r=1$  and $m=2$.

Before we present main results of the linear complexity and $k$-error linear complexity for $(f(u))_{u\ge 0}$, we prove some auxiliary statements. Define
$$ D_l(X)= \sum\limits_{u\in D_l}X^u \in \mathbb{F}_2[X]$$
for $0\leq l < p$.

\begin{lemma}\label{map}
For $r\ge 1$, $0\le l<p$ and $1\le j\le r$, the map $u\mapsto u \bmod p^j$ from $D_l$ to $\mathbb{Z}_{p^j}^*$ is surjective  and each element in $\mathbb{Z}_{p^j}^*$ exactly has $p^{r-j}$ many pre-images in $D_l$.
\end{lemma}
Proof. For each $1\le v<p^r$ with $\gcd(v,p)=1$, the numbers $v+mp^r$ belong to different $D_l~ (0\leq l < p)$ when $m$ runs through the set $\{0,1,\ldots,p-1\}$ by  Theorem \ref{HH}, hence each $D_l$ is of the form
$$
D_l=\{v+m_{lv}p^r : 1\le v<p^r, \gcd(v,p)=1, m_{lv}=v(H_{r-1}(v)-l)\bmod p \}.
$$
We will find that
$$
D_l \bmod p^{r}=\{u \bmod p^r : u\in D_l\}=\mathbb{Z}_{p^r}^*, ~~0\leq l < p,
$$
furtherly we have
$$
D_l \bmod p^{j}=\mathbb{Z}_{p^r}^* \bmod p^{j} =\mathbb{Z}_{p^j}^*, ~~0\leq l < p,
$$
for $1\le j\le r-1$. So the map $u\mapsto u \bmod p^j$ from $D_l$ to $\mathbb{Z}_{p^j}^*$ is surjective and the number of pre-images of each element in $\mathbb{Z}_{p^j}^*$ can be calculated easily. ~\hfill $\square$

From the proof of Lemma \ref{map},  each $D_l$ has the cardinality $|D_l|=p^r(p-1)$. Here and hereafter, we use $|S|$ to denote the cardinality of a set $S$.

\begin{lemma}\label{theta-pr}
Let $r\ge 2$ and $\theta \in \overline{\mathbb{F}}_{2}$ with $\theta^{p^{r}}=1$ but $\theta^{p}\neq 1$.  For $0\leq l < p$, we have
$$
D_l(\theta)=0.
$$
\end{lemma}
Proof.  We have
$$
\sum\limits_{u=0}^{p^r-1}\theta^{u}=\frac{1-\theta^{p^r}}{1-\theta}=0
$$
and
$$
\sum\limits_{u=0}^{p^{r-1}-1}\theta^{pu}=\frac{1-\theta^{p^r}}{1-\theta^p}=0.
$$
Then by Lemma \ref{map}, we derive
$$
D_l(\theta)=\sum\limits_{u\in\mathbb{Z}_{p^r}^*}\theta^{u}=\sum\limits_{u=0}^{p^r-1}\theta^{u}-
\sum\limits_{u=0}^{p^{r-1}-1}\theta^{pu}=0.
$$
We complete the proof. ~\hfill $\square$

\begin{lemma}\label{theta-p}
Let $r\ge 2$ and $\theta \in \overline{\mathbb{F}}_{2}$  with $\theta^{p}=1$. For $0\leq l < p$, we have
$$
D_l(\theta)=\left\{
\begin{array}{ll}
0, & \mathrm{if}~ \theta=1,\\
1, & \mathrm{otherwise}.
\end{array}
\right.
$$
\end{lemma}
Proof. For $\theta\neq 1$, using Lemma \ref{map} with $j=1$  we have
$$
 D_l(\theta)=p^{r-1}\sum\limits_{u=1}^{p-1}\theta^{u}=\sum\limits_{u=0}^{p-1}\theta^{u}-\theta^0
 =\frac{1-\theta^{p}}{1-\theta}+1=1.
 $$
For $\theta= 1$, we have $D_l(1)=p^r(p-1)=0$ since $|D_l|=p^r(p-1)$.  ~\hfill $\square$

\begin{lemma}\label{poly-p}
Let $\theta \in \overline{\mathbb{F}}_{2}$ with $\theta^{p}=1$ but $\theta\neq 1$ and $G(X)\in \F_2[X]$ with $1\le \deg(G(X))<p$. If  $2$ is a
primitive root modulo $p$, we have
$$
G(\theta)=1  \Longleftrightarrow G(X)=X+X^2+\ldots+X^{p-1}.
$$
\end{lemma}
Proof. Since $2$ is a
primitive root modulo $p$, we see that $1+X+X^2+\ldots+X^{p-1}$ is the minimal irreducible polynomial with the root $\theta$.
So if $G(\theta)=1$, we derive
$$
(1+X+X^2+\ldots+X^{p-1}) | (G(X)-1).
$$
With the restriction on $\deg(G(X))$, we get $G(X)=X+X^2+\ldots+X^{p-1}$. The converse is true after simple calculations.
~\hfill $\square$

Now we present our main results. We only assume $r\ge 2$ here since we have considered the case of $r=1$ in \cite{CN}, where we have more general results.

\begin{theorem}\label{klc-2-primitive}
Let $r\ge 2$ and $(f(u))_{u\ge 0}$ be the binary sequence of period $p^{r+1}$ defined in (\ref{f-binary}) using the highest-level sequence of Euler quotients in (\ref{def}) and (\ref{zero}) and a non-empty subset $\mathcal{I}$ of $\{0,1,\ldots,p-1\}$ with $1\le |\mathcal{I}|\le (p-1)/2$. If $2$ is a primitive root modulo $p^2$, then
the $k$-error linear complexity  over $\F_2$ of $(f(u))_{u\ge 0}$  satisfies
\[
 LC^{\F_2}_k((f(u))_{u\ge 0})=\left\{
\begin{array}{ll}
p^{r+1}-p^r+p-1, & \mathrm{if}\,\ 0\le k<p^{r-1}, \\
p^{r+1}-p^r+1, & \mathrm{if}\,\  p^{r-1}\le k<p^{r-1}(p-1), \\
p^{r+1}-p^r,   &\mathrm{if} ~~ p^{r-1}(p-1)\le k<p^{r-1}(p-1)|\mathcal{I}|, ~ |\mathcal{I}|>1,\\
0, & \mathrm{if}\,\ k\ge (p-1)|\mathcal{I}|, \\
\end{array}
\right.\\
\]
if $|\mathcal{I}|$ is odd, and otherwise
\[
 LC^{\F_2}_k((f(u))_{u\ge 0})=\left\{
\begin{array}{ll}
p^{r+1}-p^r, & \mathrm{if}\,\ 0\le k<p^{r-1}(p-1)|\mathcal{I}|, \\
0, & \mathrm{if}\,\ k\ge p^{r-1}(p-1)|\mathcal{I}|.
\end{array}
\right.
\]
\end{theorem}
Proof.  Let
\begin{equation}\label{Fk}
F_k(X)=\sum\limits_{l\in \mathcal{I}}D_{l}(X)+e(X)\in \F_2[X]
\end{equation}
be the generating polynomial of the sequence obtained from $(f(u))_{u\ge 0}$ by changing exactly $k$ terms of $(f(u))_{u\ge 0}$ per period,
where $e(X)$ is the corresponding error polynomial with $k$ terms. $F_0(X)$ is in fact the generating polynomial of  $(f(u))_{u\ge 0}$.
It is easy to see that if $k$ equals to or larger than the Hamming weight of $(f(u))_{u\ge 0}$, the error linear complexity will reduce to zero. So we always suppose that $k<p^{r-1}(p-1)|\mathcal{I}|$ due to $|D_l|=p^r(p-1)$, in this case $F_k(X)$ is non-zero.
We will consider the common roots of $F_k(X)$ and $X^{p^{r+1}}-1$, the number of the common roots will help us to derive the values of $k$-error linear
complexity of $(f(u))_{u\ge 0}$  by  (\ref{licom}).

We divide all roots of $X^{p^{r+1}}-1$ into four groups
$$
\mathcal{G}_1=\{\theta \in \overline{\mathbb{F}}_{2} : \theta^{p^{r+1}}=1, \theta^{p^r}\neq 1\},~~
\mathcal{G}_2=\{\theta \in \overline{\mathbb{F}}_{2} : \theta^{p^{r}}=1, \theta^{p}\neq 1\},
$$
$$
\mathcal{G}_3=\{\theta \in \overline{\mathbb{F}}_{2} : \theta^{p}=1, \theta\neq 1\},~~
\mathcal{G}_4=\{ 1\}.
$$
It is easy to check that $|\mathcal{G}_1|=p^{r+1}-p^r$, $|\mathcal{G}_2|=p^{r}-p$ and $|\mathcal{G}_3|=p-1$.

First, all $\theta\in \mathcal{G}_1$ are roots of $\Phi(X)=1+X^{p^{r}}+X^{2p^{r}}+\ldots+X^{(p-1)p^{r}}$,
which is irreducible since $2$ is a primitive root modulo $p^2$. If $F_k(\theta)=0$ for some $\theta\in \mathcal{G}_1$, we have
$$
\Phi(X)|F_k(X)
$$
and write
\begin{equation}\label{pi}
F_k(X)\equiv \Phi(X)\pi(X) \pmod {X^{p^{r+1}}-1}.
\end{equation}
Using the fact that
$$
X^{p^{r}}\Phi(X)  \equiv \Phi(X) \pmod {X^{p^{r+1}}-1},
$$
we restrict $\deg(\pi(X))<p^{r}$ and write
$$
\pi(X)=X^{v_0}+X^{v_1}+\ldots+X^{v_{t-1}} ~ \mathrm{with} ~ 0\le v_0<v_1<\ldots<v_{t-1}<p^{r},
$$
where $t\ge 1$ since $F_k(X)$ is a nonzero polynomial. Then the exponent of each monomial in $\Phi(X)\pi(X) \bmod {X^{p^{r+1}}-1}$ forms the set
$$
\{v_j+lp^r : 0\le j\le t-1, 0\le l\le p-1\},
$$
which can be divided into two sets $A$ and $B$ with
$$
A=\{v_j+lp^r : 0\le j\le t-1, 0\le l\le p-1, v_j\neq 0, q_{p,w}(v_j+lp)\in \mathcal{I}\},
$$
$$
B=\{v_j+lp^r : 0\le j\le t-1, 0\le l\le p-1\}\setminus A.
$$
By Theorem \ref{HH},  $A$ contains $|A|$ many numbers with
$$
|A|=\left\{
\begin{array}{cl}
(t-1)|\mathcal{I}|, & \mathrm{if}~ v_0=0, \\
t|\mathcal{I}|, & \mathrm{otherwise},
\end{array}
\right.
$$
and $B$ contains $tp-|A|$ many numbers.

Hence, from (\ref{Fk}) and (\ref{pi}), we find that
the set of the exponents of monomials in $e(X)$ is
$$
( \cup_{l\in \mathcal{I}}D_l\setminus A)\cup B,
$$
the cardinality of which is
$$
p^{r-1}(p-1)|\mathcal{I}|-|A|+|B|=p^{r-1}(p-1)|\mathcal{I}|+tp-\left\{
\begin{array}{cl}
2(t-1)|\mathcal{I}|, & \mathrm{if}~ v_0=0, \\
2t|\mathcal{I}|, & \mathrm{otherwise}.
\end{array}
\right.
$$
Due to $|\mathcal{I}|\le (p-1)/2$ and $tp-2t|\mathcal{I}|>0$ we have
$$
p^{r-1}(p-1)|\mathcal{I}|-|A|+|B|>p^{r-1}(p-1)|\mathcal{I}|>k,
$$
However, it is impossible that $e(X)$ has $p^{r-1}(p-1)|\mathcal{I}|$ many terms and $k$ terms simultaneously,
a contradiction. So $\Phi(X)\nmid F_k(X)$, i.e.,
\begin{equation}\label{G1}
F_k(\theta)\neq 0, ~\mathrm{for}~ \theta\in \mathcal{G}_1.
\end{equation}

Second, we consider the case $\theta\in \mathcal{G}_2$. By Lemma \ref{theta-pr} we get
\begin{equation}\label{G2}
F_k(\theta)=\left\{
\begin{array}{cl}
0, & \mathrm{if}~ k= 0,\\
e(\theta), & \mathrm{otherwise},
\end{array}
\right. ~\mathrm{for}~ \theta\in \mathcal{G}_2.
\end{equation}

Finally, we consider the case $\theta\in \mathcal{G}_3\cup \mathcal{G}_4$. By Lemma \ref{theta-p} we get
\begin{equation}\label{G3}
F_k(\theta)=\left\{
\begin{array}{cl}
|\mathcal{I}|, & \mathrm{if}~ k= 0,\\
e(\theta)+|\mathcal{I}|, & \mathrm{otherwise},
\end{array}
\right.  ~\mathrm{for}~ \theta\in \mathcal{G}_3
\end{equation}
and
\begin{equation}\label{G4}
F_k(\theta)=\left\{
\begin{array}{cl}
0, & \mathrm{if}~ k= 0,\\
e(\theta), & \mathrm{otherwise},
\end{array}
\right. ~\mathrm{for}~ \theta\in \mathcal{G}_4.
\end{equation}

Now we conclude that

(i). If $|\mathcal{I}|$ is even, we find that
$$
LC^{\F_2}_0((f(u))_{u\ge 0})=LC^{\F_2}((f(u))_{u\ge 0})=(p-1)p^r
$$
and for any $1\le k<p^{r-1}(p-1)|\mathcal{I}|$, the  number of the common roots of $F_k(X)$ and $X^{p^{r+1}}-1$ will not increase by (\ref{G1})-(\ref{G4}). So we have
$$
LC^{\F_2}_k((f(u))_{u\ge 0})=(p-1)p^r, ~~\mathrm{for}~ k<p^{r-1}(p-1)|\mathcal{I}|.
$$

(ii). If $|\mathcal{I}|$ is odd, we find that
$$
LC^{\F_2}_0((f(u))_{u\ge 0})=LC^{\F_2}((f(u))_{u\ge 0})=(p-1)p^r+p-1.
$$
Since $2$ is a primitive root modulo $p^2$, we see that
$$
\theta, \theta^2,\theta^{2^2}, \ldots, \theta^{2^{p(p-1)}-1}\in\mathcal{G}_2
$$
are different for any $\theta\in \mathcal{G}_2$.  If $e(\theta)\neq 0$ for some $\theta\in \mathcal{G}_2$, we have $e(\theta^{2^i})\neq 0$ for $0\le i<p(p-1)$.
That is to say, if such case occurs, there will be at least $p^2-p$ many $\theta\in \mathcal{G}_2$ such that $e(\theta)\neq 0$ and hence  the  number of the common roots of $F_k(X)$ and $X^{p^{r+1}}-1$ will not increase by (\ref{G2}). So according to (\ref{G2})-(\ref{G4}), we need to find the smallest $k>0$ such that the error polynomial $e(X)$ (with $k$ terms) satisfies
\begin{equation}\label{E1}
 e(\theta)=\left\{
\begin{array}{cl}
0, & \mathrm{if} ~ \theta\in \mathcal{G}_2,\\
1, & \mathrm{if} ~ \theta\in \mathcal{G}_3,\\
1, & \mathrm{if} ~ \theta\in \mathcal{G}_4,
\end{array}
\right.
\end{equation}
and
\begin{equation}\label{E2}
e(\theta)=\left\{
\begin{array}{cl}
0, & \mathrm{if} ~ \theta\in \mathcal{G}_2,\\
1, & \mathrm{if} ~ \theta\in \mathcal{G}_3,\\
0, & \mathrm{if} ~ \theta\in \mathcal{G}_4,
\end{array}
\right.
\end{equation}
respectively.

We firstly search for $e(X)$ satisfying (\ref{E1}) and consider $e(X)$ modulo $(X^{p^{r}}-1)$.  We note that
$e(X)\not\equiv 0 \bmod (X^{p^{r}}-1)$ since $e(\theta)=1$  for $\theta\in \mathcal{G}_3\cup \mathcal{G}_4$.
Let
$$
\Lambda(X):=\frac{X^{p^r}-1}{X^p-1}=1+X^{p}+X^{2p}+\ldots+X^{(p^{r-1}-1)p}\in\F_2[X].
$$
Clearly  $\Lambda(\theta)=0$ for all $\theta\in \mathcal{G}_2$ and $\Lambda(\theta)=1$ for all $\theta\in \mathcal{G}_3\cup \mathcal{G}_4$. The facts that
$$
X^p\Lambda(X)\equiv \Lambda(X) \bmod (X^{p^{r}}-1)
$$
and
$$
e(X)\equiv \tau(X)\Lambda(X) \bmod (X^{p^{r}}-1)
$$
for some non-zero polynomial $\tau(X)$ with degree $<p$ guarantee that
 the error polynomial $e(X)$ with the smallest $k>0$ terms satisfying (\ref{E1}) should be of the form
$$
e(X)\equiv \Lambda(X) \equiv 1+X^{p}+X^{2p}+\ldots+X^{(p^{r-1}-1)p} \bmod (X^{p^{r}}-1)
$$
and hence $k=p^{r-1}$. That is, when $k=p^{r-1}$ one can choose a suitable $e(X)$ as above such that the  number of the common roots of $F_k(X)$ and $X^{p^{r+1}}-1$ is equal to $p^r-1$, and for any $k<p^{r-1}$, any $e(X)$ with $k$ terms will not satisfy (\ref{E1}), this implies
that the  number of the common roots of $F_k(X)$ and $X^{p^{r+1}}-1$ will not increase (compared to the case $k=0$). So we derive
$$
LC^{\F_2}_k((f(u))_{u\ge 0})=(p-1)p^r+p-1 ~~ \mathrm{for} ~~ k<p^{r-1}
$$
and
$$
LC^{\F_2}_k((f(u))_{u\ge 0})=(p-1)p^r+1  ~~ \mathrm{for} ~~ k=p^{r-1}.
$$

Now we consider $e(X)$ satisfying (\ref{E2}). Following a similar way above, we derive by Lemma \ref{poly-p} that
the error polynomial $e(X)$ with the smallest $k>0$ terms satisfying (\ref{E2}) should be of the form
\begin{eqnarray*}
e(X)& \equiv & (X+X^2+\ldots+X^{p-1})\Lambda(X) \\
&\equiv & (X+X^2+\ldots+X^{p-1})(1+X^{p}+X^{2p}+\ldots+X^{(p^{r-1}-1)p}) \bmod (X^{p^{r}}-1)
\end{eqnarray*}
and hence  the smallest $k=p^{r-1}(p-1)$. That is, when $k=p^{r-1}(p-1)$ a suitable $e(X)$ as of the form above guarantees that the largest number of the common roots of $F_k(X)$ and $X^{p^{r+1}}-1$ is equal to $p^r$. So we derive
 $$
LC^{\F_2}_k((f(u))_{u\ge 0})=(p-1)p^r+1   ~~ \mathrm{for} ~~p^{r-1}\le k< p^{r-1}(p-1),
$$
and
$$
LC^{\F_2}_k((f(u))_{u\ge 0})=(p-1)p^r   ~~ \mathrm{for} ~~ p^{r-1}(p-1)\le k<p^{r-1}(p-1)|\mathcal{I}|, |\mathcal{I}|>1.
$$
We complete the proof.   ~\hfill $\square$

Theorem \ref{klc-2-primitive} indicates the binary sequences are cryptographically strong. By the way, we mention here the sequences $(F^{(i)}(u))_{u\ge 0}$ of Fermat quotients of order $i\ge 1$ (\ref{highorder}) and a construction of binary sequences defined by $(F^{(i)}(u))_{u\ge 0}$. It is easy to check that
$$
F^{(i)}(v+kp^i) \equiv F^{(i)}(v)-kv^{p-2} \pmod p.
$$
Following a similar proof of Theorem \ref{LC-p}, we obtain
$$
LC^{\F_p}((F^{(i)}(u))_{u\ge 0})=p^{i}+p-1
$$
for $i\ge 1$, also see a proof in \cite{L} for a more general case.  Define
$$
\widetilde{D}^{(i)}_l=\{u: 0\le u< p^{i+1},~ \gcd(u,p)=1,~F^{(i)}(u)=l\}
$$
for $l=0,1,\ldots,p-1$ and the $p^{i+1}$-periodic binary sequence $(f^{(i)}(u))_{u\ge 0}$ by
$$
f^{(i)}(u)=\left\{
\begin{array}{ll}
1, & \mathrm{if}\,\ u\bmod {p^{i+1}}\in \cup_{l\in \mathcal{I}} \widetilde{D}^{(i)}_l,\\
0, & \mathrm{otherwise},
\end{array}
\right. \quad u\ge 0,
$$
 where $\mathcal{I}$ is a non-empty subset of $\{0,1,\ldots,p-1\}$ with $1\le |\mathcal{I}|\le (p-1)/2$. If $2$ is a primitive root modulo $p^2$, using a similar proof of Theorem \ref{klc-2-primitive} we have for $i\ge 2$
\[
 LC^{\F_2}_k((f^{(i)}(u))_{u\ge 0})=\left\{
\begin{array}{ll}
p^{i+1}-p^i+p-1, & \mathrm{if}\,\ 0\le k<p^{i-1}, \\
p^{i+1}-p^i+1, & \mathrm{if}\,\  p^{i-1}\le k<p^{i-1}(p-1), \\
p^{i+1}-p^i,   &\mathrm{if} ~~ p^{i-1}(p-1)\le k<p^{i-1}(p-1)|\mathcal{I}|, ~ |\mathcal{I}|>1,\\
0, & \mathrm{if}\,\ k\ge (p-1)|\mathcal{I}|, \\
\end{array}
\right.\\
\]
if $|\mathcal{I}|$ is odd, and otherwise
\[
 LC^{\F_2}_k((f^{(i)}(u))_{u\ge 0})=\left\{
\begin{array}{ll}
p^{i+1}-p^i, & \mathrm{if}\,\ 0\le k<p^{i-1}(p-1)|\mathcal{I}|, \\
0, & \mathrm{if}\,\ k\ge p^{i-1}(p-1)|\mathcal{I}|.
\end{array}
\right.
\]
For $i=1$ the result above also holds, see \cite{CN}.

\section{Concluding Remarks}

In this paper, we define a new quotient, which coincides with the highest-level sequence of Euler quotients
decomposed as $p$-adic numbers. We use this quotient to determine the exact values of linear complexity of  the highest-level sequence of Euler quotients and values of $k$-error linear complexity
for binary sequences derived from the highest-level sequences.

We note that there are $p^{r}(p-1)|\mathcal{I}|$ many 1's in one period of the constructed binary sequences. such sequences are not balanced. It is more frequent to define binary balanced  sequences for some special applications. Unfortunately, we can't construct balanced sequences in the way described in this paper when $r>1$. However, we can modify the definition to reduce the imbalance as much as possible by defining
 $$
\widetilde{f}(u)=\left\{
\begin{array}{ll}
1, & \mathrm{if}\,\ u\bmod {p^{r+1}}\in \cup_{l\in \mathcal{I}} D_l\cup P,\\
0, & \mathrm{otherwise},
\end{array}
\right. \quad u\ge 0,
$$
where $P=\{ip : 0\le i<p^r\}$ and a non-empty subset $\mathcal{I}$ of $\{0,1,\ldots,p-1\}$ with $1\le |\mathcal{I}|\le (p-1)/2$. Together with
$$
\sum\limits_{u=0}^{p^r-1}\theta^{up}=\left\{
\begin{array}{ll}
0, & \mathrm{if}~ \theta^{p^{r+1}}=1 ~ \mathrm{but} ~ \theta^{p}\neq 1,\\
1, & \mathrm{if}~ \theta^p=1,
\end{array}
\right. ~~\mathrm{for} ~~\theta \in \overline{\mathbb{F}}_{2},
$$
we can get exact values of $k$-error linear complexity of
 $(\widetilde{f}(u))_{u\ge 0}$ if $2$ is a primitive root modulo $p^2$  by following the same way of the proof of Theorem \ref{klc-2-primitive}.

\section*{Acknowledgements}
The authors wish to thank Arne Winterhof for sending us his student's thesis, the original version of \cite{L}.

Z. Niu was partially supported by the National Natural Science Foundation of China (grants No. 61272096 and 61202395), Shanghai Municipal Natural Science Foundation (grants No.13ZR1416100 and 12ZR1443700) and the State Scholarship Fund of China Scholarship Council. ~~
Z. Chen was partially supported by the National Natural Science
Foundation of China (grant No.61373140) and the State Scholarship Fund of China Scholarship Council. ~~ X. Du was partially supported by the National Natural Science
Foundation of China (grant No.61202395) and the Program for New Century Excellent Talents in University (NCET-12-0620).

\end{document}